\documentclass[11pt]{article}
\usepackage[top=1in,bottom=1in,left=1.5in,right=1.5in]{geometry}
\usepackage[english]{babel}
\usepackage{graphicx,epstopdf,epsfig}
\usepackage{amsfonts,epsfig,fancyhdr,graphics, hyperref,amsmath,amssymb}
\usepackage{latexsym}
\usepackage[utf8]{inputenc}
\usepackage{bibentry}
\usepackage{scalerel}
\newtheorem{theorem}{Theorem}

\newtheorem{corollary}[theorem]{Corollary}

\newtheorem{con}[theorem]{Conjecture}
\numberwithin{equation}{section} 
\newcommand{\R}{\mathbb{R}}
\newcommand{\C}{\mathbb{C}}

\newcommand{\cP}{\mathcal{P}}
\newcommand{\cF}{\mathcal{F}}
\newcommand{\m}{\mathbf{m}}
\newcommand{\cM}{\mathcal{M}}
\newcommand{\cC}{\mathcal{C}}
\date{}

\begin{document}
\bibliographystyle{plain}

\title{Proof of a conjecture on polynomials preserving nonnegative matrices}
\author{Raphael Loewy\\
Department of Mathematics\\
Technion -- Israel Institute of Technology\\
\tt{loewy@technion.ac.il}}
\maketitle

\begin{abstract}
We consider polynomials in $\R[x]$ which map the set of nonnegative matrices of a given order into itself. Let $n$ be a positive integer and define $\cP_{n}:=\{p \in {\R}[x]: p(A)\geq 0, for~  all~  A \geq 0, A \in {\R}^{n,n}\}$. This set plays a role in the Nonnegative Inverse Eigenvalue Problem. Clark and Paparella conjectured that $\cP_{n+1}$ is strictly contained in $\cP_{n}$. We prove this conjecture.
\end{abstract}



\section{Introduction}\label{intro-sec}

\bigskip
The Nonnegative Inverse Eigenvalue Problem (NIEP) asks when is a list $\Lambda=  \left( \lambda_{1}, \lambda_{2}, \dots, \lambda_{n} \right) $ of complex numbers the spectrum of an $n \times n$ nonnegative matrix $A$. If it is, $\Lambda$ is said to be realizable and $A$ is a realizing matrix for $\Lambda$. This is one of the most difficult problems in matrix analysis, and is fully solved only for $n \leq 4$.

\medskip
In order to make progress solving NIEP it is essential to obtain necessary conditions on $\Lambda$ to be realizable. We demonstrate it with one example. Let
\begin{equation*}
s_{k}(\Lambda)=\sum_{i=1}^{n} \lambda_{i}^{k},~~ k=1,2,\ldots.
\end{equation*}
For $\Lambda$ to be realizable by some $A \geq 0$ it is necessary that $s_{k}(\Lambda) \geq 0, k=1,2,\ldots$, because $s_{k}(\Lambda)=trace(A^{k})$.
Stronger necessary conditions were obtained independently by Johnson \cite{CRJ} and by Loewy and London \cite{LL}.
\begin{equation}\label{JLL}
s_{k}(\Lambda)^{m} \leq n^{m-1}s_{km}(\Lambda),~~ k,m,=1,2,\ldots.
\end{equation}
Motivated by the need to obtain additional necessary conditions Loewy and London \cite{LL} defined the following set. Given a positive integer $n$, let
\begin{equation}\label{defpn}
\cP_{n}=\{p \in {\C}[x]: p(A)\geq 0, for~  all~  A \geq 0, A \in {\R}^{n,n}\}.
\end{equation}
It follows from the \eqref{JLL} that for any $p \in \cP_{n}$ and realizable $\Lambda$
\begin{equation*}
s_{k}(p(\Lambda))^{m} \leq n^{m-1}s_{km}(p(\Lambda)),~~ k,m,=1,2,\ldots,
\end{equation*}
where $p(\Lambda):= \left( p(\lambda_{1}), p(\lambda_{2}), \ldots, p(\lambda_{n}) \right)$.

\medskip
It is obvious that if $p \in \cP_{n}$ then all its coefficients are real, so we assume it throughout. It was noted in \cite{LL} that it
is possible for a polynomial to belong to some $\cP_{n}$, while having at least one negative coefficient. Bharali and Holtz \cite{BH} considered a larger class of functions, namely entire functions. They defined
\begin{equation*}
\cF_{n}=\{f~ entire: f(A)\geq 0, for~  all~  A \geq 0, A \in {\R}^{n,n}\}.
\end{equation*}
They also considered entire functions acting on certain subclasses of the set of $n \times n$ nonnegative matrices. In particular, a full characterization of $\cF_{2}$ is obtained in \cite{BH}. Clark and Paparella \cite{CP1,CP2} obtained several results concerning $\cP_{n}$. In \cite{CP1} they obtained necessary conditions for a polynomial p
to belong to $\cP_{n}$, in terms of linear inequalities that have to be satisfied by its coefficients. It is straightforward to see that for, any positive integer $n$, $\cP_{n+1} \subseteq \cP_{n}$. In \cite{CP1} and \cite{CP2} it is shown, respectively, that $\cP_{2} \subset \cP_{1}$ and $\cP_{3} \subset \cP_{2}$, and the following conjecture is raised
in \cite{CP1}.

\begin{con}\label{CPconj}
For every positive integer $n$, $\cP_{n+1} \subset \cP_{n}$.
\end{con}
In the next section we prove this conjecture.
\section{Main result}\label{main}

\begin{theorem}\label{mainth}
Let $n \geq 2$ be a positive integer, and let $p_{a}$ be the polynomial in $\R[x]$ defined by
\begin{equation*}
p_{a}[x]=1+x+x^{2}+ \cdots +x^{n-1}-ax^{n}+x^{n+1}+x^{n+2}+ \cdots +x^{2n-1}+x^{2n}.
\end{equation*}
Then, there exists $0 < a$ such that  $p_{a} \in  \cP_{n}$.
\end{theorem}
\noindent {\bf Proof:}
Let $A=(a_{i,j})$ be any $ n \times n$ nonnegative matrix. Let $B=(b_{i,j})=A^{n}$ and $C=(c_{i,j})=p_{a}(A)$. We show first that for any $0 < a \leq 2$ all the main diagonal entries of $C$ are nonnegative. It suffices to show that $c_{1,1} \geq 0$. We have $C \geq I_{n}-aB+B^{2}$,  where $I_{n}$ is the identity matrix of order $n$. Since $B^{2}_{1,1} \geq b_{1,1}^{2}$, we conclude that $c_{1,1} \geq 1-ab_{1,1}+b_{1,1}^{2}$, which is nonnegative for any $0 < a \leq 2$.

\medskip
Next, we consider the off-diagonal entries of $C$, and by symmetry it suffices to consider $c_{1,2}$. Let $ 1 \leq j \leq n$ and $x=A_{1,2}^{j}$, the $1,2$ entry of $A^{j}$.

\medskip
Then,
\begin{equation*}
x=\sum_{i_{1},i_{2} \ldots,i_{j-1}=1}^{n} a_{1,i_{1}}a_{i_{1},i_{2}} \cdots a_{i_{j-1},2}.
\end{equation*}
It is convenient to define $i_{0}=1$ and $i_{j}=2$, and rewrite
\begin{equation}\label{defofx}
x=\sum_{i_{1},i_{2} \ldots,i_{j-1}=1}^{n} \prod_{k=1}^{j} a_{i_{k-1},i_{k}},
\end{equation}
that is, $x$ is a sum of $n^{j-1}$ terms, each being a product of $j$ entries of $A$. Call each summand of \eqref{defofx} a \textit{monomial of $A$ of length $j$}. Let $\cM_{j}$ be the set of all these monomials. For each $\m \in \cM_{j}$, \textit{we keep intact the order of its $j$ factors}.

\medskip
Suppose, for example, that $j=5$ and $n=6$. Then we consider $a_{\scaleto {1,1}{6pt}}a_{\scaleto {1,4}{6pt}}a_{\scaleto {4,6}{6pt}}a_{\scaleto {6,1}{6pt}}a_{\scaleto {1,2}{6pt}}$ and $a_{\scaleto {1,4}{6pt}}a_{\scaleto {4,6}{6pt}}a_{\scaleto {6,1}{6pt}}a_{\scaleto {1,1}{6pt}}a_{\scaleto {1,2}{6pt}}$
as two distinct monomials that appear in the expansion of $A_{1,2}^{5}$, although they describe the same real number.

\medskip
One can consider the entries of $A$ as the weights (where weight zero is allowed) of the simple (that is, no multiple edges), direct and complete
graph $DK_{n}$ on the vertices $\{1,2 \ldots,n\}$. Thus, $\m=  \prod_{k=1}^{j} a_{i_{k-1},i_{k}} \in \cM_{j}$ corresponds to the path, \textit{denoted by $P(\m)$},
$1=i_{0} \rightarrow i_{1} \rightarrow i_{2} \rightarrow i_{3} \rightarrow \cdots \rightarrow i_{j-1} \rightarrow i_{j}=2$.


\medskip
Consider now $\m \in \cM_{n}$ given by $\prod_{k=1}^{n} a_{i_{k-1},i_{k}}$, where we define $i_{0}=1$ and $i_{n}=2$. Then $P(\m)$, the directed path in $DK_{n}$ corresponding to $\m$, must contain a loop or a cycle, that is, a closed path consisting of at least two edges. In the latter case, we assume the cycle is simple (that is, does not properly contain another cycle). In terms of \textit{$P(\m)$}, this means that either there exists $ 0 \leq p \leq n-1$ such that $i_{p}=i_{p+1}$, or there exist $ 0 \leq p,q \leq n$ with $p+1 < q$ and $i_{p}=i_{q}$, but $i_{r} \neq i_{s}$ for all $p \leq r < s \leq q$ with $s-r <q-p$.

\medskip
Given $\m \in \cM_{n}$, denote by $l(\m)$ the smallest number of edges in a cycle contained in \textit{$P(\m)$}. For example, if  \textit{$P(\m)$} contains a loop then $l(\m)=1$. We must have $ 1\leq l(\m) \leq n-1$.

\medskip
For $k=1,2, \ldots,n-1$, let $\cM_{n,k}$ consist of all $\m \in \cM_{n}$ with $l(\m)=k$. Then, $\cM_{n,k}, k=1,2, \ldots, n-1$, form a partition of $\cM_{n}$.

\medskip
Let $k$ be fixed, and consider  $\cM_{n,k}$. We want to define maps $\varphi_{k}$ and $\psi_{k}$ from $\cM_{n,k}$ into $\cM_{n+k}$ and $\cM_{n-k}$, respectively. For this purpose, consider $\m= \prod_{r=1}^{n} a_{i_{r-1},i_{r}} \in \cM_{n,k}$ (where, as before, $i_{0}=1$ and $i_{n}=2$). By definition, \textit{$P(\m)$} contains no cycles of length $k-1$ or less. Among the cycles of length $k$ in \textit{$P(\m)$}, consider the one occurring first. Suppose it is given by $\cC:~i_{p} \rightarrow i_{p+1} \rightarrow i_{p+2} \rightarrow \cdots \rightarrow i_{p+k-1} \rightarrow i_{p+k}=i_{p}$. Note that $\cC$ can occur consecutively more than once. Let

\begin{equation}\label{defofz}
z=\prod_{j=1}^{k} a_{i_{p+j-1},i_{p+j}}.
\end{equation}
Then we can write
\begin{equation}\label{factorization}
\m=\prod_{r=1}^{p} a_{i_{r-1},i_{r}}~ z \prod_{s=p+k+1}^{n} a_{i_{s-1},i_{s}}.
\end{equation}
Let
\begin{equation}\label{defphipsi}
\begin{aligned}
& \varphi_{k}(\m):=\prod_{r=1}^{p} a_{i_{r-1},i_{r}}~ z^{2} \prod_{s=p+k+1}^{n} a_{i_{s-1},i_{s}} \in \cM_{n+k},\\
& \psi_{k}(\m):=\prod_{r=1}^{p} a_{i_{r-1},i_{r}} \prod_{s=p+k+1}^{n} a_{i_{s-1},i_{s}}  \in \cM_{n-k}.
\end{aligned}
\end{equation}
Thus, \textit{$P(\varphi_{k}(\m))$} (respectively, \textit{$P(\psi_{k}(\m))$}), is obtained from \textit{$P(\m)$}, by adding $\cC$ again right after its first occurrence (respectively, deleting $\cC$ where it appears first).
Note that if $\cC$ appears consecutively $q$ times in \textit{$P(\m)$}, then the factor $z$ appears consecutively
$q$ (respectively $q+1$, $q-1$) times in $\m$ (respectively $\varphi_{k}(\m)$, $\psi_{k}(\m)$).

\medskip
For example, suppose that $n=12$ and
\begin{equation*}
\m=a_{\scaleto {1,3}{6pt}}a_{\scaleto {3,5}{6pt}}a_{\scaleto {5,7}{6pt}}a_{\scaleto {7,3}{6pt}}a_{\scaleto {3,10}{6pt}}a_{\scaleto {10,6}{6pt}}a_{\scaleto {6,4}{6pt}}a_{\scaleto {4,12}{6pt}}a_{\scaleto {12,10}{6pt}}a_{\scaleto {10,7}{6pt}}a_{\scaleto {7,5}{6pt}}a_{\scaleto {5,2}{6pt}}.
\end{equation*}
Then, $l(\m)=3$ and
\begin{equation*}
\begin{aligned}
& \varphi_{3}(\m)=a_{\scaleto {1,3}{6pt}}(a_{\scaleto {3,5}{6pt}}a_{\scaleto {5,7}{6pt}}a_{\scaleto {7,3}{6pt}})^{2}a_{\scaleto {3,10}{6pt}}a_{\scaleto {10,6}{6pt}}a_{\scaleto {6,4}{6pt}}a_{\scaleto {4,12}{6pt}}a_{\scaleto {12,10}{6pt}}a_{\scaleto {10,7}{6pt}}a_{\scaleto {7,5}{6pt}}a_{\scaleto {5,2}{6pt}}, \\
& \psi_{3}(\m)=a_{\scaleto {1,3}{6pt}}a_{\scaleto {3,10}{6pt}}a_{\scaleto {10,6}{6pt}}a_{\scaleto {6,4}{6pt}}a_{\scaleto {4,12}{6pt}}a_{\scaleto {12,10}{6pt}}a_{\scaleto {10,7}{6pt}}a_{\scaleto {7,5}{6pt}}a_{\scaleto {5,2}{6pt}}.
\end{aligned}
\end{equation*}

Suppose that $\m_{1} \in \cM_{n,k}$ satisfies $  \varphi_{k}(\m)=  \varphi_{k}(\m_{1})$. Then, the first occurrence of a cycle of length $k$ in \textit{$P(\m_{1})$} must be the cycle $\cC$, and since the action of $\varphi_{k}$ is to add this cycle again and keep everything else unchanged, we must have $\m_{1}=\m$, that is, $\varphi_{k}$ is one-to-one on $\cM_{n,k}$.

\medskip
Next, we consider $\m_{2}=\psi_{k}(\m)$. By its definition, it is clear that  $\psi_{k}$ is not one-to-one. We want to give an upper bound on the number of monomials in $\cM_{n,k}$ whose image under  $\psi_{k}$ is equal to $\m_{2}$. Let
\begin{equation*}
\mu(n,k)=(n-k+1)\prod_{j=1}^{k-1}(n-j).
\end{equation*}
Then, $\mu(n,k)$ is an upper bound on the cardinality of the pre-image of $\m_{2}$ under $\psi_{k}$. An element in this pre-image is obtained by inserting one monomial in $\cM_{k}$ corresponding to a cycle of length $k$. This cycle can be inserted into $\m_{2}$ in at most $n-k+1$ places. One vertex of the cycle is uniquely determined by the position where the cycle is inserted. There are $\prod_{j=1}^{k-1}(n-j)$ different cycles satisfying this requirement.

\medskip
We claim that there exists $a >0$ such that $c_{1,2} \geq 0$. Note that the only negative summands in the expansion of $c_{1,2}$ originate in
$b_{1,2}$ (recall that $B=A^{n}$). Let $\m= \prod_{r=1}^{n} a_{i_{r-1},i_{r}} \in \cM_{n}$ (where, as before, $i_{0}=1$ and $i_{n}=2$).
There exists a unique $k$, $1 \leq k \leq n-1$, such that $\m \in \cM_{n,k}$. Let $\cC$ be the first cycle of length $k$ appearing in \textit{$P(\m)$} and is given by ${i}_{p} \rightarrow i_{p+1} \rightarrow i_{p+2} \rightarrow \cdots \rightarrow i_{p+k-1} \rightarrow i_{p+k}=i_{p}$. Suppose that $z$, $\m$, $\varphi_{k}(\m)$ and $\psi_{k}(\m)$ are given by \eqref{defofz}, \eqref{factorization} and
\eqref{defphipsi}, respectively. Let
\begin{equation*}
u=\prod_{r=1}^{p} a_{i_{r-1},i_{r}},~~v=\prod_{s=p+k+1}^{n} a_{i_{s-1},i_{s}}.
\end{equation*}
Then, $\m=uzv$, $\varphi_{k}(\m)=uz^{2}v$ and $\psi_{k}(\m)=uv$. Hence,
\begin{equation}\label{quadineq}
\frac{\psi_{k}(\m)}{\mu(n,k)}- a\m +\varphi_{k}(\m)=uv(\frac{1}{\mu(n,k)}-az+z^{2}),
\end{equation}
and this expression is nonnegative for every $a \leq \frac{2}{\sqrt {\mu(n,k)}}$, since $u$ and $v$ are nonnegative. Every summand of $b_{1,2}$ is covered in such a nonnegative sum. Since $\mu(n,k) \geq 2$ for every $ 1 \leq k \leq n-1$, it follows that $p_{a} \in \cP_{n}$ for every $0 < a \leq min_{1 \leq k \leq n-1} \frac{2}{\sqrt {\mu(n,k)}}$.
\hspace{\fill} $\Box$

\medskip
Using Theorem \ref{mainth} we can prove Conjecture \ref{CPconj}.

\begin{corollary}\label{conistrue}
Let $n$ be a positive integer such that $n \geq 2$. Then, $\cP_{n+1} \subset \cP_{n}$.
\end{corollary}
\noindent {\bf Proof:}
Let $0 <a \leq min_{1 \leq k \leq n-1} \frac{2}{\sqrt {\mu(n,k)}}$. Then $p_{a} \in \cP_{n}$ by Theorem \ref{mainth}. However, by Corollary 4.2 and Theorem 4.11 in \cite{CP1} it follows that $p_{a} \notin \cP_{n+1}$.
\hspace{\fill} $\Box$

\medskip
As indicated, the cases $n=1$ and $n=2$ in the conjecture were verified by Clark and Paparella in \cite{CP1} and \cite{CP2}.

\paragraph*{Question}
Consider again $p_{a}$. The bound on $a$ given in Theorem \ref{mainth} probably is not optimal, and it is of interest to find the best bound.

\medskip
\begin{corollary}\label{morepoly}
Let $n \geq 2$ be a positive integer, and let $f_{a}$ be the polynomial in $\R[x]$ defined by
\begin{equation*}
f_{a}[x]=d_{0}+d_{1}x+ \cdots +d_{n-1}x^{n-1}-ax^{n}+d_{n+1}x^{n+1}+ \cdots +d_{2n-1}x^{2n-1}+d_{2n}x^{2n},
\end{equation*}
where $a$ is positive and $d_{i}, 0 \leq i \leq 2n, i\neq n$, are positive.
Then, there exists $0 < a$ such that  $f_{a} \in  \cP_{n}$.
\end{corollary}
\noindent {\bf Proof:}
The proof follows along the lines of the proof of Theorem \ref{mainth}, so we give only a sketch of it.

\medskip
Let $A=(a_{i,j})$ be any $ n \times n$ nonnegative matrix, and let $B=(b_{i,j})=A^{n}$. Then,
\begin{equation*}
f(A)_{1,1} \geq d_{0}-ab_{1,1}+d_{2n}b_{1,1}^{2}.
\end{equation*}
Hence $f(A)_{1,1} \geq 0$ for every $a$ such that $0 < a \leq 2 \sqrt{d_{0}d_{2n}}$.

\medskip
Next, consider $f(A)_{1,2}$. The analogue of \eqref{quadineq} is
\begin{equation*}
uv(\frac{d_{n-k}}{\mu(n,k)}-az+d_{n+k}z^{2}),
\end{equation*}
and this expression is nonnegative for every $a \leq 2 \sqrt{\frac{d_{n-k}d_{n+k}}{\mu(n,k)}}$, since $u$ and $v$ are nonnegative.
Every summand of $b_{1,2}$ is covered in such a nonnegative sum. Define $\mu(n,n)=1$. It follows that $p_{a} \in \cP_{n}$ for every $0 < a \leq min_{1 \leq k \leq n}~ 2 \sqrt{\frac{d_{n-k}d_{n+k}}{\mu(n,k)}}$.
\hspace{\fill} $\Box$




\bibliographystyle{plain}

\end{document}